\documentclass[12pt,a4paper]{amsart}
\usepackage{bbm}
\usepackage{amscd,amssymb,amsthm,amsmath}
\usepackage{graphicx}
\usepackage{jipam}

\let\set\mathbbm

\title{A Computer Proof of Tur\'an's Inequality}

\author{Stefan Gerhold}
\address{Stefan Gerhold, RISC, J. Kepler University Linz, Austria}
\email{stefan.gerhold@risc.uni-linz.ac.at}
\author{Manuel Kauers}
\address{Manuel Kauers, RISC, J. Kepler University Linz, Austria}
\email{manuel.kauers@risc.uni-linz.ac.at}
\thanks{The authors were supported by the grant F1305 of the Austrian
  FWF (Gerhold) and the grant D/03/40515 of the German DAAD (Kauers)}


\begin{document}

\begin{abstract}
We show how Tur\'an's inequality $P_n(x)^2-P_{n-1}(x)P_{n+1}(x)\geq
0$ for Legendre polynomials and related inequalities can be proven
by means of a computer procedure. 
The use of this procedure simplifies the daily work with inequalities. 
For instance, we have found the stronger inequality
$|x|P_n(x)^2-P_{n-1}(x)P_{n+1}(x)\geq 0$, $-1\leq x\leq 1$,
effortlessly with the aid of our method. 
% We strengthen Tur\'an's inequality $P_n(x)^2-P_{n-1}(x)P_{n+1}(x)\geq 0$
% for Legendre polynomials $P_n(x)$ by showing that 
% and give new proofs for some related inequalities.
% Our main point is not so much in the results themselves, but in the
% presentation of a procedure that proves inequalities involving
% recursively defined objects in an automatic way.
\end{abstract}

\maketitle

\section{Introduction}

Tur\'an showed in a 1946 letter to Szeg\H o that
\begin{equation}\label{eq:turan}
  \Delta_n(x) := P_n(x)^2-P_{n-1}(x)P_{n+1}(x)\geq0,\qquad x\in[-1,1],\ n\geq1,
\end{equation}
where $P_n(x)$ denotes the $n$-th Legendre polynomial.
Szeg\H o~\cite{szego48} gave four non-trivial proofs.
Several authors have proven analogous statements for other families of orthogonal polynomials,
and there is now a substantial body of literature~\cite{leclerc98}
devoted to these and related results.
The aim of the present note is to describe a computer algebra proof of
Tur\'an's inequality that requires as input only the three term recurrence
of the Legendre polynomials and the first two polynomials.
Our method~\cite{gerhold05} is applicable to many other inequalities,
including the following refinement of Tur\'an's result, which appears to be new.

\begin{theorem}\label{thm:refined turan}
  Let $P_n(x)$ denote the $n$-th Legendre polynomial. Then
  \begin{equation}\label{eq:str turan}
    |x|P_n(x)^2-P_{n-1}(x)P_{n+1}(x)\geq 0, \qquad x\in [-1,1],\ n\geq 1,
  \end{equation}
  with equality holding if and only if either $x=0$ and $n$ is even, or $|x|=1$.
\end{theorem}

\section{The Proving Method}

We exemplify our proving method on the classical Tur\'an inequality~\eqref{eq:turan}.
%The other inequalities claimed in this paper as well as many other
%elementary inequalities can be proven automatically in precisely the
%same way.
The idea underlying the method is complete induction
on~$n$. That is, we establish the induction step
\begin{equation}\label{eq:ind step}
  \Delta_n(x) \geq 0 \,\Longrightarrow\, \Delta_{n+1}(x) \geq 0, \qquad x\in[-1,1], n\geq1,
\end{equation}
%
%\begin{alignat}1\label{eq:ind step}
% &P_n(x)^2-P_{n-1}(x)P_{n+1}(x)\geq0
% \,\Longrightarrow\,
% P_{n+1}(x)^2-P_{n}(x)P_{n+2}(x)\geq0
%\end{alignat}
%($x\in[-1,1]$, $n\geq1$),
%
and afterwards we verify that the original inequality holds for~$n=1$.
For proving~\eqref{eq:ind step} automatically, we construct a
so-called Tarski formula whose truth implies the validity of the
induction step. Tarski formulas are quantified formulas built via logical connectives
from polynomial equations and inequalities over the reals. 
Upon replacing
$P_{n-1}(x)$, $P_n(x)$, $P_{n+1}(x)$, $P_{n+2}(x)$ in~\eqref{eq:ind step} by indeterminates
$Y_{-1}$, $Y_0$, $Y_1$, $Y_2$, we obtain the formula
\[
 \Phi:= \bigl(
 \forall\ Y_{-1},Y_0,Y_1,Y_2\in\set R:
   Y_0^2-Y_{-1}Y_1\geq0
   \,\Longrightarrow\,
   Y_1^2-Y_0Y_2\geq0 \bigr).
\]
Three things have to be remarked about this formula.
(i)~It can be decided algorithmically whether or not a given Tarski
formula is true. The classical decision procedure of
Tarski~\cite{tarski51} as well as the more efficient method of
Cylindrical Algebraic Decomposition (CAD)
due to Collins~\cite{collins75} are available for this purpose.
(ii)~If $\Phi$ holds, then~\eqref{eq:ind step} is also true, for if the implication
$Y_0^2-Y_{-1}Y_1\geq0\,\Longrightarrow\,Y_1^2-Y_0Y_2\geq0$ 
holds for all real numbers, then it holds in particular for any real
number $P_{n+i}(x)$ ($n$~and $x$ arbitrary) in place of~$Y_i$.
(iii)~Of course, $\Phi$ is false. 

In order to make the proof go through, additional knowledge about the
Legendre polynomials has to be encoded into the hypothesis part of
formula~$\Phi$.
Remarkably enough, in case of Tur\'an's inequality it is sufficient to
throw in the inequality's domain of validity and the classic 
recurrence~\cite{szego75}
\[
  (n+2)P_{n+2}(x) = (2n+3)x P_{n+1}(x) - (n+1)P_n(x), \qquad n\geq 0,
\]
of the Legendre polynomials. This requires additional
indeterminates~$N$ (representing~$n$) and $X$ (representing~$x$). 
The refined formula is
\begin{alignat*}1
 \forall\ N, X, Y_{-1},&Y_0,Y_1, Y_2\in\set R:
  \bigl(N\geq1\land -1 \leq X\leq 1\land{}\\
  &
  (N+2)Y_2=-(N+1)Y_0 + (3X+2NX)Y_1\land\\&
  (N+1)Y_1=-NY_{-1} + (X +2NX)Y_0
  \bigr)\\\,\Longrightarrow\,&
  \bigl(
   Y_0^2-Y_{-1}Y_1\geq0
   \,\Longrightarrow\,
   Y_1^2-Y_0Y_2\geq0
  \bigr).
\end{alignat*}
Using CAD, this formula can be easily verified by the computer, and by the
remarks above we may regard this as a computer proof for the fact that
Tur\'an's inequality holds for $n+1$ whenever it holds for~$n$.

To complete the proof, we have to consider the induction base $n=1$. 
Since $P_0(x)=1$, $P_1(x)=x$, and $P_2(x)=(3x^2-1)/2$, we just have to verify
the obvious formula
\[
 \forall\ X\in\set R: -1 \leq X \leq 1 \,\Longrightarrow\,
  \tfrac12(1-X^2)\geq0,
\]
which we can again leave to the computer, if we want.

Strict positivity of $\Delta_n(x)$ for $-1<x<1$ can be shown
analogously.

\section{Remarks and Further Applications}

We have to dispel any hopes that our method yields a decision procedure for
inequalities involving orthogonal polynomials or other special
functions.
Needless to say, there are many special functions that do not fit into
our recursive framework. Roughly speaking, our procedure requires
functions of~$n$ (and possibly other real parameters) such that the $n$-th value
depends polynomially on a finite number, independent of~$n$, of previous values.
For instance, the Bernoulli polynomials~$B_n(x)$ cannot be handled, since
their recurrence ``goes all the way back''\kern-.5pt.
Even if an inequality is in the input class, our method may be doomed to failure
because the sufficient condition that we check might not be satisfied
although the conjectured inequality is true.
In some cases the user can remedy this by inputting extra
equations or inequalities that the functions in question satisfy.
A third reason for failure are excessive computing time and memory overflows;
this is what happened when we tried to reprove
Gasper's extension~\cite{gasper72} of Tur\'an's inequality to Jacobi
polynomials.
Using Mathematica's implementation of CAD, we ran out of memory~(3\,GB)
after having spent forty hours of CPU time~(1.5\,GHz).
This is in contrast to the computation time needed for Tur\'an's
original inequality, whose proof was completed in just a second. The
reason for this discrepancy are the additional two parameters
appearing in the Jacobi polynomials. 

In view of the doubly exponential complexity of CAD, it is surprising
that our method is able to verify quite a few inequalities from the
literature with a reasonable amount of time.
For instance, it is a matter of seconds to verify Tur\'an's
inequality also for the following quantities in place of~$P_n(x)$:
\begin{itemize}
\item Hermite polynomials $H_n(x)$ (for $x\in\set R$),
\item Laguerre polynomials $L^\alpha_n(x)$ (for $x>0,\alpha>0$),
\item normalized Laguerre polynomials $L^\alpha_n(x)/L^\alpha_n(0)$
  (for $x\geq0,\alpha>-1$),
\item differentiated Legendre polynomials $P_n'(x)$ (For $-1\leq x\leq 1$; the inequality actually holds for all
  $x\in\set R$, but our method fails outside $[-1,1]$).
\end{itemize}
None of these results are new. We can also prove the inequality
\[
  \Delta_n(x) \geq \frac{n-1}{n+1} \Delta_{n-1}(x), \qquad x\in[-1,1],n\geq 2,
\]
which is due to Constantinescu~\cite{constantinescu05}.

Our method lends itself to playing around with conjectured inequalities;
this is how Theorem~\ref{thm:refined turan}
was obtained. Note that the absolute value function can be easily accommodated
by Tarski formulas. The cases where we claim equality in~\eqref{eq:str turan}
follow from the well-known facts $P_n(1)=1$, $P_n(-1)=(-1)^n$,
and
\[
  P_n(0) =
  \begin{cases}
    0, & n\ \text{odd}, \\
    \frac{(-1)^{n/2}}{2^n} \binom{n}{n/2}, & n\ \text{even}.
  \end{cases}
\]
These can also be proven automatically by the method described
above. However, there are a lot of established methods available for which
proving identities like these is offendingly trivial
\cite{salvy94,mallinger96,kauers03a}. 

We believe that our method could become a helpful tool for researchers
working with inequalities. 
It might not be capable of proving difficult inequalities that are of
interest in their own right (Tur\'an's inequality seems to be
exceptional in this respect), but it might be helpful for proving
elementary inequalities that appear as subproblems in the proof of
more involved theorems.

%TODO: General remarks: what to do if this does not work (extend
%indhyp and/or add extra knowledge).
%
%For $C^m_n(x)/C^m_n(1)$ and
%$P^{(\alpha,\beta)}_n(x)/P^{(\alpha,\beta)}_n(1)$ we run out of
%memory (MMA Crashes after 35h, several attempts).
%
%For $I_n(x), J_n(x)$ (Bessel) the method does not seem to terminate.
%
%For $B_n(x), E_n(x)$ (Bernoulli, Euler) the method is not applicable
%because these sequences fail to obey suitable recurrences that could
%be used.

\bibliographystyle{siam}
\bibliography{turan3}

\end{document}